%% file: deform.tex
\documentclass{amsart}
\include{headerstemplate}

\begin{document}
\title{Deformations of covers, Brill-Noether theory, and wild ramification}
\author{Brian Osserman}
\begin{abstract}
In this paper, we give a simple description of the deformations of a map
between two smooth curves with partially prescribed branching, in the cases
that both curves are fixed, and that the source is allowed to vary. Both
descriptions work equally well in the tame or wild case. We then apply this
result to obtain a positive-characteristic Brill-Noether-type result for 
ramified maps from general curves to the projective line, which even holds 
for wild ramification indices. Lastly, in the special case of rational
functions on the projective line, we examine what we can say as a result 
about families of wildly ramified maps.
\end{abstract}
\thanks{This paper was partially supported by fellowships from the
National Science Foundation and the Japan Society for the Promotion of
Science.}
\maketitle

\section{Introduction}

In studying ramified maps of curves, questions frequently arise which demand
fixing ramification on the source, or branching on the target. In the case
that the target curve is $\P^1$, the former is treated by the theory of
linear series, which naturally works up to automorphism of the image, so
we will refer to this as the linear series perspective. In contrast, we will
refer to working with fixed branching on the target (and typically allowing 
the source curve itself to vary) as the branched covers perspective. 
Often, and particularly in the context of degeneration arguments over $\C$
\cite{e-h1}, these perspectives have been considered more or less
interchangeable. However, recently a number of fundamental differences have
come to light (see for instance \cite[Prop. 5.4, Rem. 8.3]{os7}), 
particularly in positive characteristic, and it has also  
proven fruitful to pass between the perspectives to take advantages of the
distinct features of each, perhaps most notably in Tamagawa's \cite{ta1}.
In this note, we examine the deformation theory of covers with partial
branching specified, and then apply it to the perspective of linear series
to obtain a ramified Brill-Noether theorem in positive characteristic, in 
the case of one-dimensional target.

Our deformation theory result is straightforward to obtain from 
extremely well-known results, but does not appear to be stated in the
literature. It is the following.

\begin{thm}\label{main-def} Given $C, D$ smooth curves over a field $k$, and
$f:C \rightarrow D$ of degree $d$, together with $k$-valued 
points $P_1, \dots, P_n$ of $C$ such that $f$ is ramified to order at 
least $e_i$ at each $P_i$ for some $e_i$, then
the space of first-order infinitesmal deformations of $f$ together with the
$P_i$ over $k$, such that $f(P_i)$ remains fixed and the $P_i$ remain ramified
to order at least $e_i$, is parametrized by
\begin{equation}\label{main-eq} H^0(C, f^* T_D (-\sum_i
(e_i-\delta_i)P_i))\oplus k^{\sum_i(1-\delta_i)}, \end{equation}
where $\delta_i=0$ if $p | e_i$ or $f$ is ramified to order higher than
$e_i$ at $P_i$ and is $1$ otherwise.

Furthermore, the space of first-order infinitesmal deformations of $C,$ the
$P_i$ and $f$, fixing $f(P_i)$ and preserving the ramification condition at
each $P_i$, is parametrized by 
\begin{equation}\label{main-eq2} \H^1(C, T_C(-\sum_i P_i) \rightarrow
f^* T_D(-\sum_i e_i P_i)) \cong k^{d(2-2g_D)-(2-2g_C)-\sum_i (e_i-1)},
\end{equation}
where the last isomorphism requires also that $f$ be separable.

Finally, both statements still hold when some $e_i$ are allowed to be $0$,
which we interpret to put no condition on the $P_i$ at all.
\end{thm}

In the case that $f$ is tame and ramified to order exactly $e_i$ at the
$P_i$, this is well-known. The main observation of the theorem, particularly
in the first case, is that in order to obtain a useful theory, it is 
important to also consider the moduli of the ramification points, even when 
the branch points remain fixed. In the classical setting, this issue does 
not arise. 

Next, if one considers the situation with $D=\P^1$, and ramification points
specified on $C$, the perspective changes from branched covers to linear 
series with prescribed ramification. From this
point of view, classical Brill-Noether theory gives a lower bound on the
dimension of the space of maps as the $P_i$ (and even $C$) are allowed to
move. The deformation theory of Theorem \ref{main-def} gives the necessary
upper bound, and allows us to conclude the following Brill-Noether result,
generalizing \cite[Thm. 4.5]{e-h1} to positive characteristic in the case 
$r=1$.

\begin{thm}\label{brill-noether} Fix $d, n$ and $e_1, \dots, e_n$, together
with $n$ general points $P_i$ on a general curve $C$ of genus $g$. Then the 
space of separable maps of degree $d$ from $C$ to $\P^1$, ramified to 
order at least $e_i$ at $P_i$, and taken modulo automorphism of the image 
space, is pure of dimension $2d-2-g-\sum_i (e_i-1)$.
\end{thm}

Finally, we note that although from the perspective of branched covers, tame
ramification is always well-behaved and wild ramification seems more
pathological, the situation is not the same from the perspective of ramified
linear series. Indeed, from this perspective a simple dimension count
justifies the fact that wildly-ramified maps always come in infinite
families. On the other hand, we have examples from \cite[Prop. 5.4]{os7} of 
cases where
tame ramification could only produce infinite families of separable maps
with fixed ramification, and as a result of Brill-Noether theorem, can exist
only for special configurations of $P_i$. We make some elementary 
observations that, at least in certain cases when $C=D=\P^1$, wildly
ramified linear series
are in fact rather well-behaved. One such result is the following.

\begin{thm}\label{wild} Fix $d,n,m$, together with $n$ general points 
$P_i$ on $\P^1$,
and $e_1, \dots, e_n$, with $e_i$ wild for $i \leq m$, and $e_i$ tame for
$i>m$, and satisfying $2d-2 = m+\sum_i (e_i-1)$. Then the dimension of the
space of separable maps of degree $d$ from $\P^1$ to $\P^1$, ramified to
order exactly $e_i$ at $P_i$ and unramified elsewhere, and taken modulo
automorphism of the image space, is exactly $m$. Moreover, if $m=1$,
$e_1=p$, and $e_i<p$ for $i>1$, this space is non-empty if and only if the
corresponding space is non-empty when one replaces $e_1=p$ with $e_1=p-1$,
and considers maps of degree $d-1$.
\end{thm}

Note that except when explicitly stated otherwise, we make no assertions
about the non-emptyness of the space of maps with given ramification.
However, the last statement of Theorem \ref{wild} certainly produces
cases of wild ramification in which for general $P_i$, the space of maps is
non-empty of the expected dimension. This is thus better behavior than the
pathological tame examples mentioned above. The subject of existence and
non-existence will be taken up in \cite{os12}. Although the proof given here
of Theorem \ref{brill-noether} is not necessary for \cite{os12}, it does
provide the only purely algebraic, instrinsically characteristic-$p$
argument for the existence and non-existence results in question.

\section*{Acknowledgements}

I would like to thank Johan de Jong for his tireless and invaluable
guidance, and Ravi Vakil for his helpful conversations. 

\section{Deformation Theory}

Let $C, D$ be smooth curves over a field $k$, and $f:C \rightarrow D$ a
morphism of degree $d>0$. We begin by reviewing some standard deformation
theory, so that we can use formal local analysis to obtain Theorem 
\ref{main-def}. It is well-known
(see, e.g., \cite[Appendix]{va2}) that the first-order infinitesmal
deformations of $f$ are parametrized by $H^0(C, f^*(T_D))$, deformations
of a pointed curve $(C,\{P_i\})$ by $H^1(C,T_C(-\sum_i P_i))$, and deformations
of $(C,\{P_i\},f)$ by the hypercohomology group 
$\H^1(C,T_C(-\sum_i P_i)\rightarrow f^*T_D)$. By the same
token, deformations of $k$-valued points $P_i$ are parametrized simply by
$H^0(\Spec k, f^* (T_C)) \cong k$ (note that this is different from the case
of deforming $C$ along with the $P_i$ because in this case, there are no
automorphisms of $C$ to mod out by). These may be verified directly
on the Cech cocycle level using the facts that $T_C$ is the sheaf of 
infinitesmal automorphisms of $C$, and that deformations of smooth, pointed 
curves are always locally trivial. In the case of pointed curves, one 
trivializes the deformation (including of the points) locally on $C$, and 
obtains the $1$-cocycle by considering the resulting transition functions, 
taking values in infinitesmal automorphisms; in order to give 
well-defined deformations of the $P_i$, these transition functions must 
vanish along them. 

To prove Theorem \ref{main-def}, we
therefore simply need to determine the locus inside $H^0(C, f^*(T_D)) \oplus
k^n$ corresponding to maps which fix $f(P_i)$ and preserve the ramification
at $P_i$, and similarly for $\H^1(C,T_C(-\sum_i P_i)\rightarrow f^*T_D)$. 
We accomplish this easily by formal local analysis.

\begin{proof}[Proof of Theorem \ref{main-def}] Since the first statement we 
are trying to prove gives a self-contained, purely local description of a
subspace of $H^0(C, f^*(T_D)) \oplus k^n$, and fixing $f(P_i)$ and the
ramification at $P_i$ are likewise purely local conditions, it suffices to
check agreement formally locally around each $P_i$ and $f(P_i)$.
Accordingly, let $s,t$ be formal coordinates at $P_i$, $f(P_i)$
respectively. We then have $f(s)=\sum_{j \geq 0} a_j s^j$ for some $a_j \in
k$ with $a_j=0$ for $j<e_i$. First, a deformation of $f$ will be of the form
$\tilde{f}(s)=f(s)+\epsilon \sum_{j \geq 0} b_j s^j$, with the vanishing
order of the $b_j$ being the vanishing order of the section of $H^0(C,
f^*(T_D))$ inducing the deformation. Since $P_i$ corresponds to $s=0$, a
deformation of $P_i$ can be written simply as $\epsilon x$ for $x \in k$. 

If we fix both $P_i$ and $f(P_i)$, then requiring that ramification of order
$e_i$ be preserved is simply equivalent to requiring that $b_j=0$ for $j <
e_i$. If we fix $f(P_i)$, but allow $P_i$ to move, the condition that
$f(P_i)$ is fixed is simply $\tilde{f}(\epsilon x)=0$, while the condition
that $f$ remain ramified to order at least $e_i$ at $P_i$ may be expressed
vanishing to order at least $e_i$ of $\tilde{f}$ when expanded around
$s-\epsilon x$. Taylor expansion yields $\tilde{f}(s)= \sum_{j\geq 0} (a_j +
\epsilon (b_j+(j+1)a_{j+1} x))(s-\epsilon x)^j$. We see that for the first
$e_i-1$ terms to vanish, we need $b_j=0$. For the $e_i$th term, which is
$j=e_i-1$, we need $b_{e_i-1}+(e_i)a_{e_i} x=0$. The condition that
$\tilde{f}(\epsilon x)=0$ may be written $b_0+a_1 x=0$. In the case $e_i>1$,
this is automatically satisfied, while for $e_i=1$, it is redundant with the
ramification condition. Now, if $e_i a_{e_i}$ is non-zero, we see that
$b_{e_i-1}$ may be chosen arbitrarily, and uniquely determines $x$, giving
the classical case of the theorem, where $\delta=1$. On the other hand, if
$e_i a_{e_i}=0$, then we must have $b_{e_i-1}=0$, but $x$ can be arbitrary,
giving the $\delta=0$ case and completing the proof of the theorem.

Given the deformation-theoretic machinery we have already recalled, the
second statement of the theorem is even easier. Indeed, in our formal-local
trivialization, we assume by construction that the sections also correspond
to the trivial deformation, so that we are in the case above that we have
fixed both $P_i$ and $f(P_i)$, where we already noted we simply find that
$b_j=0$ for all $j<e_i$, which is equivalent to saying that our $0$-cochain
of $f^*T_D$ must vanish to order $e_i$ at $P_i$, as desired. If $f$ is 
separable, the map $T_C(-\sum_i P_i) \rightarrow f^* T_D(-\sum_i e_i P_i))$
is injective, with cokernel equal to the skyscraper sheaf of length 
$\delta_{P_i}-(e_i-1)$ at each $P_i$, where $\delta_{P_i}$ is the order of
the different of $f$ at $P_i$. The Riemann-Hurwitz formula then gives the
desired value for the dimension of the deformation space.

Finally, it is also clear from our construction that both statements work 
when some $e_i=0$ and no condition is placed on the corresponding $P_i$.
\end{proof}

\section{Brill-Noether Theory}

In this section, we prove Theorem \ref{brill-noether}, using classical
Brill-Noether theory together with Theorem \ref{main-def}.
Given two smooth curves $C,D$ over a finite-type $k$-scheme $S$, and integers 
$d,n$ and $e_1,
\dots e_n$, we have a moduli scheme $MR:=MR(C,D,d,e_1,\dots e_n)$
parametrizing $(n+1)$-tuples $(f, P_1, \dots P_n)$, where $f:C \rightarrow
D$ is a separable morphism of degree $d$, the $P_i$ are distinct sections
of $C$, and $f$ is ramified to order at least $e_i$ at $P_i$ (and
possibly elsewhere); see \cite[Appendix]{os7}. This comes with natural
forgetful morphisms $\ram: MR \rightarrow C^n$ and $\branch: MR \rightarrow
D^n$ giving the portion of the ramification and branch loci of the map $f$
which is mandated by the definition of $MR$; that is to say, the $P_i$ and
$f(P_i)$ respectively. If we fix a scheme-valued point $b$ in $D^n$,
$\branch^{-1}(b)$ is then the locus of maps $f:C \rightarrow D$ with the
specified branching above each of the $n$ points corresponding to $b$.

Recall that in
the case that $D=\P^1$, $MR$ admits a natural quotient scheme
$\overline{MR}$ which parametrizes the appropriate linear series on $C$,
together with ramification sections;
that is to say, $\overline{MR}$ represents the quotient functor obtained
simply by modding out by postcomposition with $\Aut(\P^1)$. This may be
realized as a classical relative $G^r_d$ scheme over the base $C^n$, with
prescribed ramification at the corresponding $n$ sections. Since this group
action fixes $\ram$, we have that $\ram$ factors through $\overline{MR}$.

We present the $g=0$ case of Theorem \ref{brill-noether} first, as a simpler
and more direct illustration of the general idea, using
the first statement of Theorem \ref{main-def} directly.
We thus specialize to the case that $C=\P^1$. It is not hard to see
that separable maps from $\P^1$ to itself of degree $d$ are parametrized by
an open subscheme of the Grassmannian $\G(1,d)$, and a ramification
condition of order $e_i$ corresponds to a Schubert cycle of codimension
$e_i-1$; see, e.g., \cite[\S 2]{os7}. Moreover, this description works in the
relative setting, so we conclude that $\overline{MR}$ has codimension
$\sum_i (e_i-1)$ in the trivial $\G(1,d)$ bundle over $C^n=(\P^1)^n$. With
these observations, we may easily prove the theorem.

\begin{proof}[Proof of Theorem \ref{brill-noether}, $g=0$ case] Since 
$\G(1,d)$ is smooth
of relative dimension $2d-2$ over $(\P^1)^n$, $\overline{MR}$ has dimension
at least $n+2d-2-\sum_i (e_i-1)$, and it follows that $\dim MR = \dim
\overline{MR} + \dim \Aut(\P^1) \geq n+2d+1 -\sum_i(e_i-1)$. On the other
hand, by Theorem \ref{main-def} if we are given an $f \in MR$ the dimension
of the tangent space of its fiber of $\branch$ is $h^0(\P^1, f^*
T_{\P^1} (-\sum_i (e_i-\delta_i)))+\sum_i(1-\delta_i)$ where $\delta_i=0$ if
$p | e_i$ or $f$ is ramified to order higher than $e_i$ at $P_i$ and is $1$
otherwise. Since $T_{\P^1} \cong \O(2)$, we have $f^* T_{\P^1} (-\sum_i
(e_i-\delta_i)) \cong \O(2d- \sum_i(e_i-\delta_i))$. If
$2d-\sum_i(e_i-\delta_i)$ is negative, Riemann-Hurwitz implies that $MR$ is
empty, and otherwise we find that our $h^0$ is given by 
$2d+1-\sum_i(e_i-\delta_i)$,
and the dimension of our tangent space by $2d+1-\sum_i(e_i-1)$. Thus $MR$
has dimension at most $n+2d+1-\sum_i (e_i-1)$, and this must give its
dimension precisely. The theorem then follows.
\end{proof}

We now consider the higher-genus case, assuming initially that $g \geq 2$. 
Instead of working over $\Spec k$, we let our base $S$ be a scheme \'etale
over the moduli space $\M_{g,0}$ over $k$, and let $C$ be the corresponding 
universal curve over $S$. Even in this relative setting, if we twist by
a sufficiently ample divisor $D$ on $C$ (since $g \geq 2$, we could use
high powers of the anticanonical sheaf), and then impose vanishing along 
$D$, we can realize 
$\overline{MR}$ as a closed subscheme of a Grassmannian bundle over 
$\Pic_S(C)\times _S C^n$,
with the ramification conditions corresponding to relative Schubert cycles.
This construction is carried out in the more general setting of limit 
linear series on families of curves of compact type by Eisenbud and Harris 
in the proof of \cite[Thm. 3.3]{e-h1}.

\begin{proof}[Proof of Theorem \ref{brill-noether}, $g>0$ case] The above
classical description again gives a dimensional lower bound for 
$\overline{MR}$,
this time as $(n+3g-3)+(2d-2-g)-(\sum_i(e_i-1))$, giving that 
$\dim MR \geq n+2g+2d-2-\sum_i(e_i-1)$. On the other hand, we see that the
tangent space to a fiber of $MR$ over a point $(\P^1)^n$ (under the branch
morphism composed with $S \rightarrow \Spec k$) is precisely a deformation 
of the corresponding curve $C_0$ with marked points $P_i$, together with the 
map to $\P^1$, fixing the $f(P_i)$ and the ramification conditions at the 
$P_i$. Since $g \geq 2$, there are no infinitesmal automorphisms to mod out
by in the corresponding deformation theory problem of Theorem 
\ref{main-def}, so we find that the tangent space of the fiber is described 
by that theorem, and thus has dimension $2d-2+2g-\sum_i (e_i-1)$. 
We find as before that the dimension of $MR$ is at most, hence exactly,
$n+2d-2+2g-\sum_i(e_i-1)$, and a fiber of $\overline{MR}$ over a general point
of $C^n$ must have dimension 
$n+2d-2+2g-\sum_i(e_i-1)-(n+3g-3)-3=2d-2-g-\sum_i(e_i-1)$, as desired.

We conclude with the $g=1$ case (we could handle the $g=0$ case similarly,
but since we have already given a proof in that case, we will not do so).
Here, we argue as when $g \geq 2$, but let $S$ be \'etale over $\M_{1,1}$. 
The Brill-Noether lower bound works as before to give the relative dimension
of $2d+1-g-\sum_i(e_i-1)$ for $MR$ over $C^n$, with the ample divisor 
in the construction arising from our given section. Because $S$ has
dimension $3g-3+1$ in this case, we find we need the fiber of $MR$ over
a point $(\P^1)^n$ to have dimension $1$ greater than before. The tangent
space of this fiber at a point $(C_0,\{P_i\}_i,f)$ now includes a $P_0$ with
$e_0=0$, giving the extra dimension, as desired.
\end{proof}

\begin{rem} One observes in the $g=1$ case above that even if all 
ramification is specified, the fiber dimension for fixed branching is 1. 
The reason for this is that the construction of $MR$ doesn't see the
marked point on the genus 1 curve which comes from a point of $S$, and 
still allows changing the ramification sections $P_1,\dots,P_n$ by 
automorphism of the underlying curve $C_0$.
\end{rem}

We remark that although the $g=0$ case of this theorem is extremely easy in 
characteristic $0$,
the situation is more delicate in characteristic $p$. In particular, in the
intersection of the ramification Schubert cycles frequently has an excess
intersection corresponding to inseparable maps of lower degree. Furthermore,
examples such as $x^{p+2}+tx^p+x$ give tamely ramified situations where all
non-empty fibers of $\ram$ have greater than the expected dimension; in such
situations, $\ram(MR)$ necessarily fails to dominate $(\P^1)^n$ even though
the expected dimension is non-negative. However, the argument of Theorem
\ref{brill-noether} implies that this can never happen for $\branch(MR)$.

\begin{cor} In the situation of Theorem \ref{brill-noether}, but 
allowing the $P_i$ and $C$ to vary as in the proof, if $f\in MR$ is any 
point with any neighborhood $U \subset MR$, then $U$ dominates 
$(\P^1)^n$ under the $\branch$ morphism, with all fibers smooth of 
dimension $2d-2+2g-\sum(e_i-1)+\epsilon_g$, where $\epsilon_g$ is the
dimension of the infinitesmal automorphism space of a curve of genus $g$ 
(i.e., 3 if $g=0$, 1 if $g=1$, and 0 otherwise).
\end{cor}

\begin{proof} We saw in the proof of Theorem \ref{brill-noether} that if
$MR$ is non-empty, it is pure of dimension 
$n+2g+2d-2+\epsilon_g-\sum_i(e_i-1)$, and that
the tangent space at any point in any fiber of the $\branch$ morphism has
dimension precisely $2d-2+2g+\epsilon_g-\sum_i(e_i-1)$. The corollary follows.
\end{proof}

\section{Wild Ramification}

We conclude with some largely elementary remarks on wild ramification, again
in the situation that $C=D=\P^1$. The first observation is that by
Riemann-Hurwitz in characteristic $p$, if any $e_i$ are wild, in order to
have separable maps with the desired ramification, we must have
$2d-2>\sum_i(e_i-1)$. The codimension count of the previous section then
implies that the separable maps with at least the specified ramification
will necessarily form an infinite family. Thus, the fact that wildly
ramified maps come in infinite families is elementary from the perspective
of linear series. With the exception of one application of Theorem
\ref{brill-noether} to prove Theorem \ref{wild}, all our observations will
be of a completely elementary nature, but we hope they may shed some light
on the behavior of wildly ramified maps.

\begin{proof}[Proof of Theorem \ref{wild}] The primary assertion follows
from Theorem \ref{brill-noether} together with the observation that under
our hypotheses, the
locus of maps $f$ with exactly the specified ramification is open in the
locus of maps with at least the specified ramification. Indeed, having
ramification exactly $e_i$ at $P_i$ is always open, since ramification can
only decrease under deformation. On the other hand, by Riemann-Hurwitz a
deformation cannot have additional ramification away from the $P_i$, since
the different at each wild $e_i$ is necessarily at least $e_i$. 

The second assertion will be a special case of the following proposition.
\end{proof}

Considering Theorem \ref{brill-noether} and the preceding argument, one
might be led to expect that the space of wildly ramified maps having higher
different at the wild points would have higher dimension. However, this is
not necessarily the case. In the argument for openness above, we use
minimality of the different in a key way. Indeed, if one deforms a map with
greater than minimal discriminant, a new ramification point specializing to
the wild point can appear, as illustrated (indirectly) by the following
two propositions.

\begin{prop}\label{foo} Let $d,n$ and $e_1, \dots, e_n$ be positive integers, with the
$e_i$ less than $p$, and $2d-2=\sum_i (e_i-1)$. Also, let $P_1, \dots, P_n$
be distinct points on $\P^1$. Then there exists a separable map of degree
$d$ from $\P^1$ to itself, ramified to order $e_i$ at $P_i$, if and only if
there exists a separable map of degree $d+p-e_1$, ramified to order $e_i$ at
$P_i$ for $i>1$, and order $p$ at $P_1$. The dimension of the space wild maps 
in this situation is $1$ more than the dimension of the space of tame ones.
\end{prop}

\begin{proof} We may assume that $P_1 = \infty$, and
$f(\infty)=\infty$. Then we can go back and forth between the wild and tame
cases simply by adding appropriate multiples of $x^p$, since the fact that 
$e_i <p$ for $i>1$ implies that the ramification away from $\infty$ will 
remain unchanged. We also use that the different in the wild case at
$\infty$ is less than $2p$, so that if we subtract a multiple of $x^p$ from
a wild map, the degree of the numerator cannot drop below the degree of the
denominator. The difference in dimension comes from the fact that the
multiple of $x^p$ added to obtain a wild map can be arbitrary.
\end{proof}

\begin{prop} Let $d,n$ and $e_1, \dots, e_n$ be positive integers, with
$e_1=d=p$, $e_i$ less than $p$ for $i>1$, and $2d-2>\sum_i (e_i-1)$. Also,
let $P_1, \dots, P_n$ be distinct points on $\P^1$. Then the space of
separable maps of degree $d$ from $\P^1$ to itself, ramified to order $e_i$
at $P_i$ and unramified elsewhere, and taken modulo automorphism of the
image, is non-empty of dimension $1$.
\end{prop}

\begin{proof} Without loss of generality, we can assume $P_1=\infty$ and
$f(\infty)=\infty$; then $f$ is given by a polynomial of degree $p$. Since
$e_i<p$ for $i>1$, the ramification conditions for $i>1$ determine the
derivative of $f$ (up to scaling). On the other hand, since
$\sum_{i>1}(e_i-1)<p-1$, an $f$ with the desired derivative always exists.
The space is $1$-dimensional because the $x^p$ term may be scaled
independently from the lower-order terms.
\end{proof}

Up until now, all of our examples have suggested that the dimension of a
wildly ramified family will always be equal to the number of wildly 
ramified points. However, the following example shows that this is not always
the case, even for families existing for general $P_i$. 

\begin{ex} The family $\frac{x^{2p}+t_1x^{p+1}+t_2}{x^p+t_1 x}$ for
$t_1,t_2$ non-zero is a
two-dimensional family of rational functions (modulo automorphism of the 
image) ramified to order $p$ at infinity, and unramified elsewhere.
\end{ex}

\begin{rem} One can try to say more about the case with a single wildly 
ramified point by generalizing the argument of Proposition \ref{foo}, 
inductively inverting as necessary and subtracting off inseparable 
polynomials. However, there are subtleties to be aware for this sort of
argument. In particular, neither the tame ramification indices nor the 
dimension
of the tame family obtained in this process will be determined by the
ramification indices and degree of the wildly ramified map. Indeed, the maps
$\frac{x^5(x^{10}+x^7-2x)+1}{x^{10}+x^7-2x}$ and
$\frac{x^5(x^5(x^5+x^4-x^3+2x)+x^2+2x+1)+x^5+x^4-x^3+2x}
{x^5(x^5+x^4-x^3+2x)+x^2+2x+1}$
in characteristic $5$ are both of degree $15$, ramified to order $5$ at
infinity, and simply ramified at the $6$th roots of unity, but the tame
functions they reduce to are $x^7-2x$ and $\frac{x^2+2x+1}{x^5+x^4-x^3+2x}$
respectively; the former moves in a one-dimensional family, while the latter
doesn't. Nonetheless, it is interesting to note that each of the wild maps 
moves in a $2$-dimensional family.

However, even in the tame case the situation of one index being at least $p$
while the others are less than $p$ is pathological, so it is not clear how
much general intuition one should attempt to draw from this case. That said,
it is interesting that at least for this example, it seems the dimension in 
the wild case is in fact
more uniform than the dimension in the tame case. This suggests that an
approach other than reducing to the tame case is likeliest to be productive 
for analyzing the dimension of families of wildly ramified maps.
\end{rem}

\bibliographystyle{hamsplain}
\bibliography{hgen}
\end{document}

%% file: headerstemplate.tex
\usepackage{amssymb,euscript,amsmath, mathrsfs}
\usepackage[dvips]{graphicx}
\usepackage[dvips]{color}

\newcounter{ENUM}

\input xy
\xyoption{all}
\CompileMatrices

\def\P{{\mathbb P}}

\def\C{{\mathbb C}}
\def\H{{\mathbb H}}
\def\G{{\mathbb G}}

\def\O{{\mathcal O}}

\def\M{{\mathcal M}}

\def\Aut{\operatorname{Aut}}

\def\ram{\operatorname{ram}}
\def\branch{\operatorname{branch}}

\def\Spec{\operatorname{Spec}}
\def\Pic{\operatorname{Pic}}

\newtheorem{thm}{Theorem}[section]
\newtheorem{prop}[thm]{Proposition}

\newtheorem{cor}[thm]{Corollary}

\theoremstyle{definition}

\newtheorem{ex}[thm]{Example}

\theoremstyle{remark}

\newtheorem{rem}[thm]{Remark}

\numberwithin{equation}{section}
\numberwithin{figure}{section}